\documentclass{amsart}

\newtheorem{prop}{Proposition}[section]  
\newtheorem{cor}{Corollary}[section]  
\newtheorem{lemma}{Lemma}[section]

\newcommand{\NN}{{\mathbb N}}
\newcommand{\ZZ}{\mathbb Z}

\def\M{{\mathcal M}}

\def\R{{\mathcal R}}

\begin{document}

\pagestyle{plain}
  
\title  
{
Corrigendum and addendum to `Linearly recurrent subshifts have a finite
number of non-periodic factors'
}  

\author{Fabien  Durand}
  
\address{Laboratoire Ami\'enois
de Math\'ematiques Fondamentales  et
Appliqu\'ees, CNRS-UMR 6140, Universit\'{e} de Picardie
Jules Verne, 33 rue Saint Leu, 80039 Amiens Cedex 01, France}
\email{fabien.durand@u-picardie.fr}

\maketitle

\begin{abstract}
We prove that a subshift $(X,T)$ is linearly recurrent if and only if it
 is a primitive and proper $S$-adic subshift. This corrects Proposition
 6 in F. Durand ({\it Ergod. Th. \& Dynam. Sys. {\bf 20}} (2000), 1061--1078).
\end{abstract}

\section{Introduction and definitions}

In this paper we freely use the definitions and the notations of \cite{Du}.
Proposition 6 in \cite{Du} is false: There exist primitive $S$-adic
subshifts that are not linearly recurrent (LR). We will give an example. Nevertheless the
other part of this proposition is true: If a subshift is LR then it is
primitive $S$-adic.

The author apologizes for the mistake. We correct this
proposition with the next one, but before we have to give and to recall some
definitions. Let $A$ and $B$ be two finite alphabets.

Let $x$ be an element of $A^{\NN}$ or $A^{\ZZ}$. 
We call {\it
occurrence} of $u\in A^*$ in $x$ every integer $i$ such that
$x_{[i,i+|u|-1]} = x_i x_{i+1} \dots  x_{i + |u| - 1}= u$.
A {\it return word} to $u\in A^*$ in $x$ is a word $w$ such
that $wu$ has an occurrence in $x$, $u$ is a prefix of $wu$ and $u$ has
exactly 2 occurrences in $wu$. We say that $x$ is {\it linearly
recurrent} (LR) (with constant $K\in\NN$) if it is uniformly recurrent
and if for all $u$ having an occurrence in $x$ and all return words, $w$,
to $u$ in $x$ we have $|w| \leq K|u|$.

Let $T$ be the shift transformation defined on $A^{\ZZ}$. 
We say the subshift $(X,T)$ is generated by $x$ if $X$ is the set of the sequences $z$ such that $z_{[i,j]}$ has an occurrence in $x$ for all intervals $[i,j]\subset \ZZ$.
The subshift $(X,T)$ is {\rm linearly recurrent} if it is minimal and contains a LR sequence. We remark that if $x\in A^{\ZZ}$ is linearly recurrent then $x$ and $x_{[0,+\infty )}$ generate the same LR subshift.

Let $a$ be a letter of $A$, $S$ a finite set of
morphisms $\sigma$ from $A(\sigma)\subset A$ to $A^{*}$ and $(\sigma_{n}
: A_{n+1} \rightarrow A_{n}^{*}; n\in \NN )$ be a sequence of $S^{\NN}$
such that $(\sigma_{0}\sigma_{1}\cdots\sigma_{n} (aa\cdots ) ; n\in\NN)$
converges in $A^{\NN}$ to $x$. We will say that $x$ is a {\it $S$-adic
sequence on $A$}  (generated by $(\sigma_i ; i\in \NN) \in S^{\NN}$ and
$a$). If there exists an integer $s_0$ such that for all non-negative
integers $r$ and all $b\in A_{r}$ and $c\in A_{r+s_0+1}$, the letter $b$
has an occurrence in $\sigma_{r+1}\sigma_{r+2}\cdots \sigma_{r+s_0}
(c)$, then we say that $x$ is a {\it primitive} $S$-adic sequence (with
constant $s_0$).

Let $\sigma : A \rightarrow B^*$ be a
morphism. We say $\sigma $ is proper if there exist two letters $r,l\in
B$ such that for all $a\in A$ the first letter of $\sigma (a)$ is $l$
and the last letter of $\sigma (a)$ is $r$. We say the sequence $x\in A^{\NN}$ is a
{\it proper $S$-adic sequence} if it is a
$S$-adic sequence and the morphisms in $S$ are proper. The subshift
generated by a proper $S$-adic sequence is called {\it proper $S$-adic subshift}.

\begin{prop}
\label{lrsadic}
The subshift $(X,T)$ is LR if and only if it is a primitive and proper
 $S$-adic subshift. 
\end{prop}

\section{Counterexample to the Proposition 6 of \cite{Du}}

In this section we give a counterexample to Proposition 6 in \cite{Du},
i.e., a primitive $S$-adic subshifts that is not LR.

Let $A = \{ a,b,c\}$ be an alphabet, and, $\sigma : A \rightarrow A^{*}$
and $\tau : A \rightarrow A^{*}$ be two morphisms defined by 
$$
\begin{array}{lllll}
\sigma (a) & = acb \ , & & \tau (a) & = abc \ , \\
\sigma (b) & = bab \ , & & \tau (b) & = acb \ , \\
\sigma (c) & = cbc \ , & & \tau (c) & = aac \ .
\end{array}
$$
We call $\M$ the set of all finite composition of elements of $S = \{\sigma
, \tau \}$. For each element $\rho $ of $\M$ there exists a unique $n\in \NN$ such
that $|\rho (a)| = |\rho (b)| = |\rho (c)| = 3^n$. We set $|\rho |
= 3^n$.
We first give a lemma which proof is left to the reader. 

\begin{lemma}
\label{second}
Let $z\in A^{\NN}$ and $n\in \NN$. The difference between two successive
 occurrences of the word $ca$ in $\sigma^n \tau (z)$ is greater than
 $3^{n+1}$.
\end{lemma}

Let $x\in A^{\NN}$ be the primitive $S$-adic sequence defined by 
$$
x= \lim_{n\rightarrow +\infty} \sigma \tau \sigma^2 \tau \cdots \sigma^n \tau (aaa \dots ). 
$$
We show $x$ is not LR. Let $n$ be an integer, we
set $\rho_n = \sigma \tau \sigma^2 \tau \cdots \sigma^n \tau $ and $
y= \lim_{l\rightarrow +\infty} \sigma^{n+2} \tau \sigma^{n+3} \tau
\cdots \sigma^{n+l} \tau (aaa \dots )$. We have $x = \rho_n
\sigma^{n+1} \tau (y)$. 

Since $\sigma (a)$ is a prefix of $y$, $\tau \sigma (a) = abcaacacb$,
$\sigma^{n+1} (a) = au$  and $\sigma^{n+1} (c) = vc$, for some $u,v \in
A^*$, then $ca$ appears in $\sigma^{n+1} \tau (y)$.
Let $w$ be a return word to $ca$ in $\sigma^{n+1} \tau (y)$. Hence $wca$
($= ca w' ca $ for some $w'$) appears in $x$. The word $\rho_n (ca)$ appears exactly twice in
$\rho_n (caw'ca)$ (the proof of this fact is left to the reader) hence $\rho_n (w)$ is a return word to $\rho_n (ca)$
in $x$.

Moreover from Lemma \ref{second} we have $|w| \geq 3^{n+2}$. It implies $x$ is not LR because
$$
\frac{|\rho_n (w)|}{|\rho_n (ca)|} 
=
\frac{|w| |\rho_n| }{2|\rho_n |}
\geq 
\frac{3^{n+2}}{2} .
$$ 

Let $z\in B^{\NN}$ where $B$ is a finite alphabet. We denote by $L (z)$ the set of all words having an occurrence in $z$. For all $n\in \NN$ we
define $p_z (n)$ to be the number of distinct words of length $n$ in
$L(z)$. In \cite{Du} it is proved that if $z$ is LR with constant $K$ then
$p_z (n) \leq K n $ for all $n\in \NN$.
Even if $x$ is not LR, there exists a constant $C$ such that $p_n (x)
\leq Cn$ for all $n\in \NN $. This is a consequence of the following proposition.

\begin{prop}
Let $A$ be a finite alphabet, $a$ be a letter of $A$ and $(\sigma_n : A_{n+1} \rightarrow A_n^+ ; n\in \NN)$ be a sequence of morphisms such that $A_n \subset A$ for all $n\in \NN$, $a\in \cap_{n\in \NN} A_n$ and 
$$
y = \lim_{n\rightarrow +\infty} \sigma_0 \sigma_1 \cdots \sigma_n (aaa\dots) .
$$
Suppose moreover $\inf_{c\in A_{n+1}} |\sigma_0 \sigma_1 \cdots \sigma_n (c)|$ tends to $+\infty$ and there exists a constant $D$ such that 
$$
|\sigma_0 \sigma_1 \cdots \sigma_{n+1} (b)|
\leq
D|\sigma_0 \sigma_1 \cdots \sigma_{n} (c)|
$$
for all $b\in A_{n+2}$ and $c \in A_{n+1}$, and all $n\in \NN$. Then $p_y (n) \leq D{\rm (Card} (A))^2 n$.
\end{prop}

\begin{proof}
This proof follows the lines of the proof of Proposition V.19 in
\cite{Qu}.

Let $n\geq 1$. The sequence $(\inf_{c\in A_{k+1}} |\sigma_0 \sigma_1
\cdots \sigma_k (c)| )_{k\in \NN}$ is non-decreasing and tends to $+\infty$, hence there exists
$p\in \NN$ such that 
$$
\inf_{c\in A_{p}} |\sigma_0 \sigma_1
\cdots \sigma_{p-1} (c)|
\leq n \leq
\inf_{c\in A_{p+1}} |\sigma_0 \sigma_1
\cdots \sigma_p (c)|.
$$
From that, every word $w\in L(y)$ of length $n$ has an occurrence $i$ in some
$\sigma_0 \sigma_1 \cdots \sigma_p (bc)$, where $b$ and $c$ are two
letters of $A$, with $i\leq |\sigma_0 \sigma_1 \cdots \sigma_p
(b)|-1$. Consequently 
$$
p_y (n) 
\leq 
({\rm Card} (A))^2 \sup_{c\in A_{p+1}} |\sigma_0 \sigma_1
\cdots \sigma_p (c)|
$$
$$
\leq
({\rm Card} (A))^2 D \inf_{c\in A_{p}} |\sigma_0 \sigma_1
\cdots \sigma_{p-1} (c)|
\leq
D ({\rm Card} (A))^2 n .
$$
This ends the proof.
\end{proof}

\medbreak

\begin{cor}
Let $A$ be a finite alphabet, $a$ be a letter of $A$, $l$ be a positive
 integer and $(\sigma_n : A \rightarrow A^* ; n\in \NN)$ be a sequence
 of morphisms of constant length $l$ and 
$$
y = \lim_{n\rightarrow +\infty} \sigma_0 \sigma_1 \cdots \sigma_n (aaa\dots) .
$$
Then $p_y (n) \leq l{\rm (Card} (A))^2 n$.
\end{cor}

\section{A sufficient condition for a primitive $S$-adic sequence to be
LR}
 
This sufficient condition is given in the following lemma and will
be used in the sequel.

\begin{lemma}
\label{condsadiq}
Let $S$ be a finite set of morphisms. Let $x$ be a primitive
 $S$-adic sequence generated by $(\sigma_i : A_{i+1} \longrightarrow A_i^* ; i\in \NN)$ and $a$ (with constant $s_0$). For all $n\in \NN$ suppose 
$
\lim_{l\rightarrow +\infty} \sigma_n \sigma_{n+1} \cdots \sigma_l (aaa\dots) 
$
exists and call it $x^{(n)}$. Let $D_n$ be the largest difference
 between two consecutive occurrences of a word of length 2 in
 $x^{(n)}$. 

If $(D_n ; n \in \NN )$ is bounded then $x$ is LR.
\end{lemma}
\begin{proof}
Let $x = \lim_{n\rightarrow +\infty} \sigma_0\sigma_1\cdots \sigma_n
(aaa\dots)$. It follows from Lemma 7 of \cite{Du} that $x$ is uniformly
recurrent. We set $S_k = \sigma_{0}\cdots \sigma_{k}$ for all $k\in
\NN$.  Let $u$ be a non-empty word of $L (x)$ such that $|u|\geq \max \{
|S_{s_0}(b)| ; b\in A_{s_0+1}\}$, and $v$ be a return word to $u$. We
denote by $k_0$ the smallest positive integer $k$ such that $|u| < \min
\{ |S_k (b)| ; b\in A_{k+1}\}$; we remark that $k_0\geq s_0+1$. There
exists a word of length 2, $bc$, of $L (x^{(k_0+1)})$ such that $u$ has
an occurrence in $S_{k_0} (bc)$. The largest difference between two
successive occurrences of $bc$ in $x^{(k_0+1)}$ is bounded by $D=
\max_{n\in \NN} D_n$ (which does not depend on $k_0$), hence we have
$$
|v|
\leq D\max \{|S_{k_0} (d)| ; d\in A_{k_0+1} \} 
\leq DK \min \{|S_{k_0} (d)| ; d\in A_{k_0+1} \}
$$
$$
\leq DK \max \{|S_{k_0 -1} (d)| ; d\in A_{k_0} \} \min \{|\sigma_{k_0} (d)| ; d\in A_{k_0+1} \}
$$
$$
\leq DK^2 \min \{|S_{k_0 -1} (d)| ; d\in A_{k_0} \} \min \{|\sigma_{k_0} (d)| ; d\in A_{k_0+1} \}
$$
$$
\leq DK^2 \min \{|\sigma_{k_0} (d)| ; d\in A_{k_0+1} \} |u|,
$$
where $K$ is the constant given by Lemma 8 of \cite{Du}, i.e., $K$ is
such that for all integers $r,s$ with $s-r\geq s_0+1$ and all $b,c$ of $A_{s+1}$ we have
$
|\sigma_{r}\cdots \sigma_{s} (b)| \leq
K |\sigma_{r}\cdots \sigma_{s}(c)|
$.
We set $M = DK^2 \max \{|\sigma_{i} (d)| ; i\in \NN , d\in A_{i+1} \} $. For all $u$ of $L (x)$ greater than $\max \{ |S_{s_0}(b)| ; b\in A_{s_0+1}\}$ and all $v$ in $\R_{u}$ we have $|v|\leq M|u|$. Hence $x$ is LR.
\end{proof}

\medbreak

\section{A necessary and sufficient condition to be LR}

In the original proof of Proposition 6 in \cite{Du} we use the notion of
return word. In the proof of Proposition \ref{lrsadic} we will do the same but we will use an extension
of this notion which was defined in \cite{DHS}. We will take a sequence $x$ belonging to $X$ and, using
these ``new'' return words, we will show that $x^+ = x_{[0,+\infty )}$ is a
primitive and proper $S$-adic sequence. The subshift $(X,T)$ being minimal we
see that it is generated by $x^+$ and, consequently, $(X,T)$ is a
primitive and proper $S$-adic subshift.

Let $A$ be a finite
alphabet, $x\in A^{\ZZ} $, and, $u$ and $v$ two words of $A^*$. We say $w \in A^*$ is a {\it return word
to $u.v$} in $x$ if there exist two consecutive occurrences
$j,k$ of $uv$ in $x$ such that $w=x_{[j+|u|,k+|u|)}$. It is immediate to
check that a word $w\in A^+$ is a return word to $u.v$ in $x$ if and
only if:

\medskip

1) $uwv$ has an occurrence in $x$, and

2) $v$ is a prefix of $wv$ and $u$ is a suffix of $uw$, and

3) the word $uwv$ contains exactly two occurrences of the word $uv$.

\medskip

We denote by $\R_{x,u.v}$
the set of return words to $u.v$ in $x$. If $u$ is the
empty word $\epsilon$ then the return words to $u.v$ are the return words to $v$ defined
in \cite{Du} and we set $\R_{x,u.v} = \R_{x,v}$. The return words to
$u.v$ are different from the return words to $\epsilon.uv$ but we have $
\# \R_{x,u.v}
= \# \R_{x,\epsilon . uv} = \# \R_{x,uv}$. 

We suppose now that $x$ is a uniformly
recurrent sequence. 
It is easy to see that for all $u,v\in L (x)$ the set
$\R_{x,u.v}$ is finite. It will be convenient to label the return
words. We enumerate the elements $w$ of $\R_{x,u.v}$ in the order of the
first appearance of $uwv$ in $x_{[-|u|,+\infty )}$. This defines a
bijective map $\Theta_{x,u.v}: R_{x,u.v}\to\R_{x,u.v} \subset A^+$ where
$
R_{x,u.v} = \{ 1,\ldots,\# \R_{x,u.v} \} 
$:
$u\Theta_{x,u.v} (k) v$ is the $k$-th word of the type $uwv$, $w\in
\R_{x,u.v}$, appearing in $x_{[-|u|,+\infty )}$.

We consider $R_{x,u.v}$ as an alphabet. The map $\Theta_{x,u.v}$ defines
a morphism from $R_{x,u.v}$ to $A^*$ and the set $\Theta_{x,u.v}
(R_{x,u.v} ^*)$ consists of all concatenations of return words to $u.v$.

The following proposition is important in the proof of Proposition \ref{lrsadic}.

\begin{prop}[\cite{DHS}]
\label{onetoone}
The map $\Theta_{x,u.v} : R_{x,u.v}^+\to A^+$ is one to one.
\end{prop}

\begin{proof}[Proof of Proposition \ref{lrsadic}]
Let $S$ be a finite set of proper morphisms and suppose $(X,T)$ is a
primitive $S$-adic subshift generated by $(\sigma_i : A_{i+1}
\rightarrow A_i^* ; i\in \NN)\in S^{\NN}$ and $a$ (with constant
$s_0$). Let 
$$
x = \lim_{n\rightarrow +\infty} \sigma_0\sigma_1\cdots \sigma_n (aaa...).
$$
We prove that $x$ is LR and consequently that the subshift it generates
is LR.

As the morphisms are proper the limit 
$
\lim_{l\rightarrow +\infty} \sigma_n \sigma_{n+1} \cdots \sigma_l (aaa\dots) 
$
exists for all $n\in \NN$. We call it $x^{(n)}$ and we define $D_n$
as in Lemma \ref{condsadiq}.

The composition of two proper morphism is again proper. Consequently,
from the primitivity, we
can suppose that $s_0 = 0$ and that for all $n\in \NN$, all $a\in
A_{n+1}$ and all $b\in A_n$ the letter $b$ appears in $\sigma_n (a)$.

Let $n\in \NN$ and set $\tau = \sigma_n \sigma_{n+1}$. It is a proper
substitution. Let $l$ and $r$ be respectively the first and the last
letter of the images of $\sigma_n$. Let $y$ be a one-sided sequence of
$A_{n+2}^{\NN}$ and $z = \sigma_n \sigma_{n+1} (y)$. The words of length 2
having an occurrence in $z$ are exactly the words of length 2 having an
occurrence in some $\sigma_n (e)$, $e\in A_{n+1}$, and the word $rl$. On
the other hand the letters of $\sigma_{n+1} (y)$ appear with gaps
bounded by $K_n = 2\max \{ | \sigma_{n+1} (e)| ; e\in A_{n+2} \}$. Consequently
the words of length 2 of $z$ appear with gaps bounded by $K_n\max \{ |
\sigma_{n} (e)| ; e\in A_{n+1} \} $ and, a fortiori, $D_n \leq K_n\max \{ |
\sigma_{n} (e)| ; e\in A_{n+1} \} $ for all $n\in \NN$. Moreover $S$ being
finite the sequence $(D_n ; n\in \NN)$ is bounded. Lemma \ref{condsadiq}
implies $(X,T)$ is LR.

We suppose now that $(X,T)$ is LR. The periodic case is trivial hence we
suppose that $(X,T)$ is not periodic. From Proposition 5 in \cite{Du} there exists $K\geq 2$ such that 
$$
(\forall u\in L (X)) (\forall w\in \R_u ) \left( |u|/K \leq |w| \leq K |u| \right),
$$
We set $\alpha = K^2(K+1)$. Let
$x= (x_n ; n\in \ZZ )$ be an element of $X$. It suffices to prove that
$x_{[0,+\infty )}$ is a primitive and proper $S$-adic sequence.

For each non-negative
integer $n$ we set $u_n = x_{-\alpha^n} \cdots x_{-2} x_{-1}$,
$v_n = x_0x_1\cdots x_{\alpha^n-1}$, $\R_n = \R_{x,u_n.v_n}$, $R_n =
R_{x,u_n.v_n}$ and $\Theta_n = \Theta_{x,u_n.v_n}$. Let $n$ be a
positive integer and $w$ be a return word to $u_n.v_n$. The word $w$ is
a concatenation of return words to $u_{n-1}.v_{n-1}$. The map
$\Theta_{n-1}$ being one to one (Proposition \ref{onetoone}), this
induces a map $\lambda_n $ from $R_n $ to $R_{n-1}^{*} $ defined by
$\Theta_{n-1}\lambda_n = \Theta_{n}$. We set $ \lambda_0 = \Theta_0
$. For each letter $b$ of $R_n$ we have $|\Theta_{n-1}\lambda_n (b)
|\leq K|u_nv_n| = 2K \alpha^n.$
Moreover each element of $\R_{n-1}$ is greater than $(2 \alpha^{n-1})/K$ hence 
$$
|\lambda_n (b) | \leq \frac{K^2\alpha^n}{\alpha^{n-1}} = \alpha K^2.
$$
By Proposition 5 of \cite{Du} we have $\# R_n = \#
\R_{x,u_nv_n} \leq K(K+1)^2$, consequently the set $M = \{ \lambda_n ;
n\in \NN \}$ is finite. The definition of $R_n$ implies that $\Theta_{n}
(1)x_0 x_1\cdots x_{\alpha^n-1}$ is a prefix of $x_{[0,+\infty )}$ for all $n\in \NN$ and 
$
\lambda_0 \lambda_1\cdots \lambda_n (1) = \Theta_{n} (1).
$
Proposition 4 of \cite{Du} implies that the length of $\Theta_{n} (1)$ tends to infinity with $n$ and
$$
x = \lim_{n\rightarrow + \infty} \lambda_0 \lambda_1 \cdots \lambda_n (11\cdots).
$$
Let $n$ be an integer greater than 1. Each word of length $2K\alpha^n$ has an occurrence in
each word of length $2K(K+1)\alpha^n$ (Proposition 5 of
\cite{Du}). Hence each element of $\R_n$ has an occurrence in each word
of length $2K(K+1)\alpha^n$. Let $w$ be an element of $\R_{n+1}$, we
have $|w| \geq 2\alpha^{n+1}/K = 2K(K+1)\alpha^n $. Therefore each element of $\R_n$ has an occurrence in each element of
$\R_{n+1}$. It means that if $b$ belongs to $R_{n+1}$ then each letter of
$R_n$ has an occurrence in $\lambda_{n+1} (b)$. Hence $x$ is a
primitive $S$-adic sequence.

It remains to show each $\lambda_n$ is proper. Let $w$ be a return word of
$\R_n$. The word $wv_{n-1}$ is a concatenation of return words to
$v_{n-1}$. Let $p\in \R_{n-1}$ be such that $pv_{n-1}$ is a prefix of
$wv_{n-1}$ and consequently of $wv_n$. We know $v_n$ is also a prefix of $wv_{n}$ and 
$$
|v_n| = \alpha^n = \alpha^{n-1} (K^3 + K^2) \geq (2K+1) \alpha^{n-1} \geq |p| + |v_{n-1}|.
$$
Consequently $pv_{n-1}$ is a prefix of $v_n$. Let $l\in R_{n-1}$ be such
that $\Theta_{n-1} (l) = p$. Then $l$ is the first letter of $\lambda_n (c)$ for
all $c\in R_n$. 

In the same way there exists $s\in \R_{n-1}$ and $r\in R_{n-1}$ such
that $\Theta_{n-1} (r) = s$ and $u_{n-1}s$ is a suffix of
$u_n$. Hence $r$ is the last letter of $\lambda_n (c)$ for
all $c\in R_n$ and $\lambda_n$ is proper. 
\end{proof}

\medbreak

\section{LR sturmian sequences}

We give a correct proof of the next proposition (which is stated in
\cite{Du}) because the original proof used Proposition 6 in \cite{Du}.

\begin{prop}
A sturmian subshift $(\Omega_{\alpha}, T)$ is LR if and only if the
 coefficients of the continued fraction of $\alpha$ are bounded.
\end{prop}
\begin{proof} 
Let $0< \alpha <1$ be an irrational real number and $[0:i_1+1 , i_2,i_3,
\dots ]$ be its continued fraction. From Proposition 9 in \cite{Du} we
know that  $\Omega_{\alpha } = \Omega (x)$ where 
$$
x =
\lim_{k\rightarrow +\infty} \tau^{i_1} \sigma^{i_2} \tau^{i_3} \sigma^{i_4} \cdots \tau^{i_{2k-1}} \sigma^{i_{2k}} (00\cdots) ,
$$
$\tau (0) = 0$,  $\tau (1) = 10$, $\sigma (0) = 01$ and $\sigma (1) = 1$.
We just have to prove
that if the coefficients of the continued fraction of $\alpha$ are
bounded then the sequence $x$ is LR. The other part of the proof is in
\cite{Du} and do not use Proposition 6 in \cite{Du}.

Let $i$, $j$ and $k$ be in $\NN^*$. We have 
$$
\tau^i \circ \sigma^j \circ \tau^k (0) = 0(10^i)^j 
\hbox{ and }
\tau^i \circ \sigma^j \circ \tau^k (1) = 10^i(0(10^i)^j)^k. 
$$
Consequently if $x$ belongs to $\{ 0,1 \}^{\NN}$ then the set of the
words of length 2 having an occurrence in $y = \tau^i \circ \sigma^j \circ
\tau^k (x)$ is $\{ 00 , 01, 10
\}$. Moreover the difference of two successive occurrences of 01
(resp. 10) in $y$ is less than $i+2$ (resp. $i+2$), and, the difference between two
successive occurrences of 00 in $y$ is less than $2j+3$ if $i=1$ and
less than 3 if $i\geq
2$. Consequently the difference between two successive occurrences of a
word of length 2 in $y$ is bounded by $2\max \{ i,j,k \} +3$. 

The same bound can be found for $\sigma^i \circ \tau^j \circ \sigma^k$. 

For all $n\in \NN$, let $x^{(n)}$ and $D_n$ be defined as in Lemma
\ref{condsadiq}. Hence, if the sequence $(i_n ; n\in \NN )$ is bounded by $K$,
then $D_n$ is bounded by $2K+3$ for all $n\in \NN$. Lemma \ref{condsadiq} ends the proof.
\end{proof}

\medbreak

{\bf Acknowledgements.}
I would like to thank Krzysztof Wargan who
pointed out the mistake in the proof of Proposition 6 in \cite{Du}.

\end{document}